\magnification=\magstep1
\def\N{I\! \! N}
\def \bull{\vrule height .9ex width .8ex depth -.1ex}
\def\Blackboard R{I \!\! R}
{\footline={$\ast$\hfil}
\hrule height 0pt
\vskip 1in
\bigskip
\bigskip
\bigskip
\centerline{NEW BOUNDS FOR HAHN AND KRAWTCHOUK POLYNOMIALS}
\smallskip
\centerline{by}
\medskip
\centerline{Holger Dette}
\centerline{Institut f\-ur Mathematische Stochastik}
\centerline{Technische Universit\"at Dresden}
\centerline{Mommsenstr. 13}
\centerline{01062 Dresden}
\centerline{GERMANY}
\medskip
\centerline{Technical Report \#93-31C}
\vskip 2in
\centerline{Department of Statistics}
\centerline{Purdue University}
\bigskip
\centerline{July 1993}   
\vfill\break
}
\pageno=1
\parskip=6pt
\baselineskip=15.8pt
\eject
\centerline{NEW BOUNDS FOR HAHN AND KRAWTCHOUK POLYNOMIALS}  
\smallskip
\centerline{by}
\smallskip
\centerline{Holger Dette$^*$}  
\centerline{Institut f\-ur Mathematische Stochastik}  
\centerline{Technische Universit\"at Dresden}
\centerline{Mommsenstr. 13}
\centerline{01062 Dresden}
\centerline{GERMANY}
\bigskip
\bigskip
\bigskip
\footnote{}{$^*$Research supported in part by the Deutsche
Forschungsgemeinschaft.}
\footnote{}{AMS Subject Classification:  33C45}
\footnote{}{Keywords and phrases:  Hahn polynomials, Hahn-Eberlein polynomials,
Krawtchouk}
\footnote{}{polynomials, dual Hahn polynomials, discrete trigonometric identity.}

\centerline{ABSTRACT}

\medskip
\noindent
For the Hahn and Krawtchouk polynomials orthogonal on the
set $\{0, \ldots,N\}$ new identities for the sum of squares are derived
which generalize the trigonometric identity for the Chebyshev
polynomials of the first and second kind.  These results are applied in
order to obtain conditions (on the degree of the polynomials)
 such that the polynomials are bounded (on the
interval $[0,N]$) by their values at the points $0$ and $N$.  As special
cases we obtain a discrete analogue of the trigonometric identity and
bounds for the discrete Chebyshev polynomials of the first and second
kind.

\vskip .7cm 

\noindent
{\bf \underbar {1. Introduction}}.  The Hahn polynomials may be defined in
terms of a hypergeometric series
$$
\eqalign{
Q_n (x,\alpha, \beta, N) & ~=~ _3F_2 
{\left(\matrix{-n,&n+\alpha+\beta+1,
&-x;&1\cr
\alpha+1, &-N\cr}\right)}\cr
& \cr
&=~\sum^n_{k=0} {{(-n)_k (n+\alpha+\beta+1)_k (-x)_k} \over {k!~
(\alpha+1)_k (-N)_k}} \ \ \ ~~~~~~(n=0,\ldots,N)\cr
}
$$
\noindent
where $\alpha, \beta > -1$, $(a)_0 = 1, (a)_k = a(a+1) \ldots (a+k-1)$.
These polynomials are limiting cases of some general systems of orthogonal
polynomials (see Hahn (1949)) and satisfy for $n, m = 0, \ldots, N$ the
orthogonality relation
$$
\sum^N_{x=0} \rho (x, \alpha, \beta, N) Q_m (x,\alpha,\beta,N) Q_n
(x,\alpha,\beta,N) ~=~ {{\delta_{nm}} \over {\pi_n (\alpha, \beta, N)}}
\leqno(1.1)
$$
\noindent
where
$$
\rho (x,\alpha, \beta, N) ~=~{\displaystyle{{x+\alpha}\choose x} 
{{N-x+\beta}\choose {N-x}} 
\over \displaystyle{{N+\alpha+\beta+1}\choose N}}
\leqno(1.2)
$$
\noindent
and
$$
\pi_n (\alpha, \beta, N) ~=~ {{(-1)^n (-N)_n (\alpha + 1)_n (\alpha +
\beta + 1)_n} \over {n!~ (N+\alpha+\beta+2)_n (\beta+1)_n}} \ \ {{2n +
\alpha + \beta + 1} \over {\alpha + \beta + 1}}.
\leqno(1.3)
$$
\noindent
For some properties and applications of the Hahn polynomials we refer the
reader to the work of Karlin McGregor (1961, 1962), Gasper (1974, 1975)
and Wilson (1970).  The polynomials $Q_n (x,\alpha, \beta, N)$ can be
seen as the discrete analogue of the Jacobi polynomials and most of the
``classical'' orthogonal polynomials can be obtained as limits from the Hahn
polynomials when the parameters tend to infinity (see Gasper (1975)).

\noindent
As an  example we  consider the Krawtchouk polynomials which can be defined
as the limit $(q=1-p, p \in (0,1))$
$$
\eqalign{
k_n (x,p,N)~ &=~ \lim_{t \rightarrow \infty} Q_n (x,pt,qt,N) ~=~ _2F_1 
(-n,-x,-N;1/p)\cr
&= ~\sum^n_{k=0} {{(-n)_k (-x)_k} \over {k!~(-N)_k }} ({1 \over p})^k\cr
}
\leqno(1.4)
$$
\noindent
and are orthogonal with respect to the jump function
$$
{N \choose x} p^x (1-p)^{N-x} \ \ \ ~~~~~~~~~x = 0, \ldots, N
\leqno(1.5)
$$
\noindent
(see Krawtchouk (1929)).
\smallskip
\noindent
In this paper we will discuss some new properties of the
orthogonal polynomials with respect to the measures (1.2) and
(1.5).  After presenting some preliminary results in Section 2 we
present new identities for squares of Krawtchouk and Hahn
polynomials in Section 3 and 4 which generalize the trigonometric
identity for the Chebyshev polynomials of the first and second
kind.  We will apply these results in order to obtain conditions 
(on the degree of the polynomials) such that
the polynomials $Q_n (x,\alpha,\beta,N)$ and $k_n(x,p,N)$ are bounded on
the interval 
$[0,N]$ by their values at the points $0$ and $N$.  For the Hahn
polynomials these bounds extend and improve 
results of Zaremba (1975) while
for the Krawtchouk polynomials it is shown that 
$$
|k_n (x,p,N)|
~\leq ~\max \ \{|k_n (0,p,N)|, |k_n (N,p,N)|\} ~=~\max \{1, ({q
\over p})^n\}   ~~~~~~~~x \in [0,N]
$$ 
whenever the degree of the polynomial satisfies $n
\leq  {N \over 2} +1$.  Similar  results are also
given for the 
 dual Hahn polynomials and
the Hahn- Eberlein polynomials.
\bigskip
\bigskip
\noindent
{\bf \underbar {2. \ Preliminaries.}}~~In this section we will briefly
discuss some general aspects of orthogonal polynomials which will
be needed in the following sections.  The notation used here is
that of Karlin and Shapely (1953) and Karlin and Studden (1966).
Let $\xi$ denote a probability measure on the interval $[0,N]$
with moments 
$$
c_j = \int^N_0 x^j d \xi (x) ~~~~~~~~~~~~~~~~~~~(j=0,\ldots,N)
$$
 and
let $P_\ell (x), Q_\ell (x), R_\ell (x), S_\ell (x)$ denote the
orthonormal polynomials with respect to the measures $d \xi (x),
x(N-x) d\xi (x), xd\xi (x)$ and $(N-x) d\xi (x)$, respectively.
The leading coefficients of these polynomials can be expressed by
ratios of the determinants
$$
\left\{\matrix{ {\underline D}_{2 \ell} (\xi) = |(c_{i+j})^\ell_{i,j=0}|&
{\overline D}_{2 \ell} = |(Nc_{i+j-1} - c_{i+j})^\ell_{i,j=1}|\cr
& \cr
{\underline D}_{2 \ell +1} (\xi) = | (c_{i+j+1})^\ell_{i,j=0}|&
{\overline D}_{2\ell+1} = |(Nc_{i+j} -
c_{i+j+1})^\ell_{i,j=0}|\cr
}\right.
\leqno(2.1)
$$
\noindent
(see e.g. Karlin and Studden (1966 p. 109)).
For a point $(c_1,\ldots,c_\ell)$ in the interior of the moment
space
\hfuzz=5pt
$$
{\cal M}_\ell = \{ (c_1,\ldots,c_\ell) | c_j = \int^N_0 x^j d \xi
(x) \hbox{ for  some   probability  measure  on } [0,N]\ ~(j=1,\ldots ,l)\}
$$
\noindent
let $(c_1,\ldots, c_{\ell - 1}, c^-_\ell)$ and $(c_1, \ldots,
c_{\ell - 1}, c^+_\ell)$ denote the boundary points of ${\cal M}_\ell$
corresponding to the lower and upper principal representation
associated with  the point $(c_1, \ldots, c_{\ell - 1}) $$\in \ int \
({\cal M}_{\ell - 1})$ (see Karlin and Studden (1966), p.55).  It is
well known  (see e.g. Karlin and Shapely (1953), p. 59)
that the quantities $c^+_\ell$ and $c^-_\ell$ can be
expressed in terms of the determinants (2.1), that is 
$$
c^+_\ell = c_\ell + {{{\overline D}_\ell (\xi)} \over {{\overline
D}_{\ell - 2} (\xi)}} \ , ~~~~~~~~\ c^-_\ell = c_\ell - {{{\underline
D}_\ell (\xi)} \over {{\underline D}_{\ell - 2} (\xi)}} \ ~~~~~~~~~\ell
\geq 1
\leqno(2.2)
$$
\noindent
where we define ${\underline D}_{-1} (\xi) = {\underline D}_0
(\xi) = {\overline D}_{-1} (\xi) = {\overline D}_0 (\xi) = 1$
(note that the ratios in (2.2) are well defined because $(c_1,
\ldots, c_{\ell - 1} ) \in \ int \ ({\cal M}_{\ell - 1})).$
\smallskip
\noindent
In the following we will make use of the determinants 
defined in (2.1)
where the moment of highest order is replaced by $c^+_{2 \ell}$
$(c^+_{2 \ell + 1})$ in the determinants ${\underline D}_{2 \ell}
(\xi)$ $ ({\underline D}_{2 \ell + 1} (\xi))$ and by $c^-_{2 \ell}$
$ (c^-_{2 \ell + 1})$ in the determinants ${\overline D}_{2 \ell}
(\xi)$ $ ({\overline D}_{2 \ell + 1} (\xi))$.  The corresponding
modified determinants are denoted by ${\underline D}^+_{2 \ell}
(\xi), {\underline D}^+_{2 \ell + 1}(\xi), {\overline D}^-_{2
\ell} (\xi)$ and ${\overline D}^-_{2 \ell + 1} (\xi)$,
respectively.  Using the representation (2.2) it is then easy to
see that
$$
\left\{\matrix{
{\underline D}^+_j (\xi) &=& {\underline D}_j (\xi)
  + \displaystyle{{\underline D}_{j-2} (\xi) \over
  {\overline D}_{j-2} (\xi)} {\overline D}_j (\xi)&
    j=2\ell,\ 2 \ell+1\cr
&&\cr
  {\overline D}_j^- (\xi) &=&
    {\overline D}_j (\xi) -
	 \displaystyle{{\overline D}_{j-2} (\xi) \over
	  {\underline D}_{j-2} (\xi)} {\underline D}_j (\xi)
	  & j= 2\ell,\ 2 \ell+1.\cr
}\right.
\leqno(2.3)
$$
\noindent
In a recent paper Dette (1994) established new identities for the
 orthonormal polynomials $P_\ell (x), Q_\ell (x),
R_\ell (x)$ and $S_\ell (x)$ with respect to the measures $d \xi
(x)$, $x (N-x) d \xi (x)$, $ x d \xi (x)$ and $(N-x) d \xi (x)$
respectively.  For example, it is shown that for any arbitrary
probability measure $\xi$ on the interval $[0,N]$ the corresponding orthonormal
polynomials satisfy the identity
$$
\eqalign{
& \sum^{n-1}_{\ell=1}
{{{\underline D}_{2 \ell - 1} (\xi)}
\over
{{\overline D}_{2 \ell - 1} (\xi)}}
\left[
{{{\overline D}_{2 \ell -2} (\xi)}
\over
{{\underline D}_{2 \ell -2} (\xi)}}
-
{{{\overline D}_{2 \ell} (\xi)}
\over
{{\underline D}_{2 \ell} (\xi)}}
\right]
P^2_\ell (x)
~+~
{{{\underline D}_{2n-1} (\xi) \ {\overline D}_{2n-2} (\xi) \
{\underline D}_{2n} (\xi)}
\over
{{\overline D}_{2n-1} (\xi) \ {\underline D}_{2n-2} (\xi) \
{\underline D}^+_{2n} (\xi)}}
P^2_n (x) \cr
\cr
&~~~~~~~~~~+ ~(N-x) \sum^{n-1}_{\ell=0}
{{{\overline D}_{2 \ell} (\xi)}
\over
{{\underline D}_{2 \ell} (\xi)}}
\left[
{{{\underline D}_{2 \ell -1} (\xi)}
\over
{{\overline D}_{2 \ell -1} (\xi)}}
- 
{{{\underline D}_{2 \ell +1} (\xi)}
\over
{{\overline D}_{2 \ell +1} (\xi)}}
\right]
S^2_\ell (x) \cr
& \cr
& ~~~~~~~~~~ 
=~
1 ~-~ x(N - x)
{{{\underline D}_{2n-1} (\xi) {\overline D}_{2n} (\xi)}
\over
{{\overline D}_{2n-1} (\xi) {\underline D}_{2n}^+ (\xi)}}
Q^2_{n-1} (x)
}
\leqno (2.4)
$$
(note that the identities were originally stated on the interval $[-1,1]$
but can easily be transferred to arbitrary intervals).
If $N=1$ and
 $$d \xi (x) = 
{{dx}
\over
{\pi \sqrt{x (1-x)}}}
$$ is the arcsin distribution, then it is
straightforward to show that ${\overline D}_{2 \ell} (\xi) =
{\underline D}_{2 \ell} (\xi) = ({1 \over 2})^{\ell (2 \ell + 1)},
{\overline D}_{2 \ell + 1} (\xi) = {\underline D}_{2 \ell + 1}
(\xi) = ({1 \over 2})^{(\ell + 1) (2 \ell + 1)}$ (see e.g. Karlin and
Studden (1966), p.123).  The
polynomials $P_\ell (x)$ and $Q_\ell (x)$ are proportional to the
Chebyshev polynomials of the first and second kind (on the
interval $[0,1]$) and the identity (2.4) reduces to the
trigonometric identity.  In this sense (2.4) can be seen as an
extension of the trigonometric identity for arbitrary orthogonal
polynomials on compact intervals.
  For the Jacobi polynomials identities of the form
(2.4) have been estiblished in Dette (1994).  In order to derive
similar results for the Hahn and Krawtchouk polynomials we need
explicit expressions for the determinants of the moment matrices
corresponding to the jump functions in (1.2) and (1.5) which will
be derived in the following sections. 
\bigskip
\bigskip

\noindent
\hfuzz=65pt
{\bf \underbar {3. \ Identities and bounds for Hahn polynomials.}}
~~It follows from (1.1) and (1.3) that the jump function in (1.2)
defines a (discrete) probability measure $\xi_\rho$ on the set $\{0,
\dots, N\}$ and that the orthonormal polynomials with respect to
the measure $d \xi_\rho (x)$ are given by $\sqrt{\pi_n (\alpha,
\beta, N)}  Q_n (x, \alpha, \beta, N)$ $ (n=0, \ldots, N)$.  Using
the elementary properties of the Gamma function and (1.1) we
obtain
$$
\eqalign{
& \sum^N_{x=0} Q_m (x-1, \alpha + 1, \beta + 1, N-2) Q_n (x-1,
\alpha + 1, \beta + 1, N-2) x (N-x) \rho (x, \alpha, \beta, N)
\cr
& \cr
&
= \sum^{N-2}_{x=0} Q_m (x, \alpha + 1, \beta + 1, N-2) Q_n (x,
\alpha + 1, \beta + 1, N-2)  \cr
&~~~~~~~~~~~~~~~~~~~~~~~~~~~~~~~~~~~~~~~~~~~~
\times\rho (x, \alpha + 1, \beta + 1,
N-2)
{{N(N-1) (\alpha + 1) (\beta + 1)}
\over
{(\alpha + \beta + 2)(\alpha+\beta+3) }}
 \cr
& \cr
&=~
{{N (N-1) (\alpha + 1) (\beta + 1)}
\over
{(\alpha + \beta + 2) (\alpha + \beta + 3)}}
\cdot
{{\delta_{m,n}}
\over
{\pi_n (\alpha + 1, \beta + 1, N-2)}}
\cr
}
$$
\noindent
which shows that the polynomials
$$
\sqrt{
{{(\alpha + \beta + 2) (\alpha + \beta + 3)}
\over 
{N(N-1) (\alpha + 1) (\beta + 1)}}
\pi_n (\alpha + 1, \beta + 1, N-2)}
~Q_n (x-1, \alpha + 1, \beta + 1, N-2) 
\leqno(3.1)
$$
\noindent
($n=0,\ldots ,N-2$) are orthonormal with respect to the measure $x(N-x) d \xi_\rho
(x)$.  Similarly it can be shown that the orthonormal polynomials
with respect to the measures $x d \xi_\rho (x)$ and $(N-x) d
\xi_\rho (x)$ are given by
$$
\sqrt{
{{(\alpha + \beta + 2)}
\over
{(\alpha + 1) N}}
\pi_n (\alpha + 1, \beta, N-1)}
~Q_n (x-1, \alpha + 1, \beta, N-1) 
\leqno(3.2)
$$
\noindent
($n=0,\ldots ,N-1$) and by
$$
\sqrt{
{{(\alpha + \beta + 2)}
\over
{(\beta + 1)N}} 
\pi_n (\alpha, \beta + 1, N-1)}
~Q_n (x, \alpha, \beta + 1, N-1) 
\leqno(3.3)
$$
\noindent
($n=0,\ldots , N-1$), respectively.
\medskip
\bigskip
\noindent
{\bf \underbar {Theorem 3.1.}}~~  For $\ell = 0, \ldots, N$ define
$h_\ell (x, \alpha, \beta, N) = 
{{(\alpha + 1)_\ell}
\over
{(\beta + 1)_\ell}} Q_\ell (x, \alpha, \beta, N)$, then the Hahn
polynomials 
satisfy the following identities

\noindent
a) For $n=0,\ldots ,N-1$~:
$$
\eqalign{
& \sum^{n-1}_{\ell=1}
{{2 \ell + \alpha + \beta + 1}
\over
{N}}
\{(\alpha + \beta + 1) (2 \ell - N) + 2 \ell^2\}
\left\{
{{(\alpha + \beta + 2)_{\ell - 1}}
\over
{(N+ \alpha + \beta + 2)_\ell}}
{N \choose \ell}
h_\ell (x, \alpha, \beta, N) \right\}^2
\cr
& \cr
&+~
\left\{ {{N-1} \choose {n-1}}
{{(\alpha + \beta + 2)_{n-1}}
\over
{(N+\alpha + \beta + 2)_{n-1}}}
h_n(x, \alpha, \beta, N) \right\}^2
\cr
 &\cr
&+~ (\beta - \alpha) (1 - {x \over N}) \sum^{n-1}_{\ell = 0}
{{2 \ell + \alpha + \beta + 2}
\over
{(\beta + 1)^2}}
\left\{
{{N-1} \choose \ell}
{{(\alpha + \beta + 2)_\ell}
\over
{(\alpha + \beta + N + 2)_\ell}}
h_\ell (x, \alpha, \beta + 1, N -1) \right\}^2
\cr
& \cr
&=~ 1 ~-~ {x \over N} (1 - {x \over N}) \left\{
{{(\alpha + \beta + 2)_n}
\over
{(\alpha + \beta + N + 2)_{n-1}}}
{{N-2} \choose {n-1}}
{
h_{n-1} (x-1, \alpha + 1, \beta + 1, N-2) \over \beta +1}
\right\}^2 \cr
}
$$
b) For $n=0,\ldots , N-1$~:
$$
\eqalign{
&\sum^n_{\ell = 1}
{{2 \ell + \alpha + \beta + 1}
\over
N}
\{ (\alpha + \beta + 1) (2 \ell - N) + 2 \ell^2 \}
\left\{ {N \choose \ell}
{{(\alpha + \beta + 2)_{\ell - 1}}
\over
{(\alpha + \beta + N + 2)_\ell}}
Q_\ell (x, \alpha, \beta, N) \right\}^2
\cr
&\cr
&+~ {x \over N} \left\{ {{N - 1} \choose n}
{{(\alpha + \beta + 2)_n}
\over
{(\alpha + \beta + N + 2)_n}}
{{\alpha + n + 1}
\over
{\alpha + 1}}
Q_n (x - 1, \alpha + 1, \beta, N-1) \right\}^2
\cr
&\cr
&+ ~(\alpha - \beta) {x \over N} \sum^{n-1}_{\ell = 0}
(2 \ell + \alpha + \beta + 2) 
\left\{ {{N-1} \choose \ell}
{{(\alpha + \beta + 2)_\ell}
\over
{(\alpha + \beta + N + 2)_\ell}}
{{Q_\ell (x-1, \alpha + 1, \beta, N-1)}
\over
{\alpha + 1}}\right\}^2
\cr
&\cr
&=~ 1 ~-~ (1 - {x \over N}) \left\{ {{N-1} \choose n}
{{(\alpha + \beta + 2)_n}
\over
{(\alpha + \beta + N + 2)_n}}
Q_n (x, \alpha, \beta + 1, N-1) \right\}^2 
\cr
}
$$
c) For $n=0,\ldots , N-2$~:
$$
\eqalign{
&x (1 - {x \over N}) \sum^{n-1}_{\ell = 0} \{ (\alpha + \beta + \ell
+ 2) (N - 2{\ell-2}) - (\ell + 1) N\}
(2 \ell + \alpha + \beta + 3) \cr
& ~~~~~~~~~~~~~~~~~~~~~~~~~~~~~~~~~~~~~~~~~~~~~~~~~~~ \times \left\{
{{h_\ell (x-1, \alpha + 1, \beta + 1, N-2)}
\over
{(\beta + 1) (N-1)}}\right\}^2
\cr
& \cr
&+~ {x \over N} (1 - {x \over N}) \left\{
{{(\alpha + \beta + 2 + n) (N-n-1)}
\over
{(\beta + 1) (N-1)}}
h_n (x-1, \alpha + 1, \beta + 1, N-2) \right\}^2
\cr
&\cr
&+~ (\beta - \alpha) (1 - {x \over N}) \sum^n_{\ell = 0} (2 \ell +
\alpha + \beta + 2) \left\{
{{h_\ell (x, \alpha, \beta + 1, N-1)}
\over
{\beta + 1}} \right\}^2
=~ 1~ -~ \{h_{n+1} (x, \alpha, \beta, N)\}^2 \ 
\cr
}
$$
d) For $n=0,\ldots , N-1$~:
$$
\eqalign{
&\sum_{\ell=1}^n {2\ell + \alpha + \beta + 1\over N} \{(\alpha +
   \beta+ 1) (2\ell - N) + 2\ell^2\} \left\{{(\alpha + \beta + 
   2)_{\ell - 1}\over (\alpha + \beta + N + 2)_\ell} {N\choose 
   \ell} h_\ell (x, \alpha, \beta, N)\right\}^2 \cr
&\cr
&+~ (\beta - \alpha) (1 - {x\over N}) \sum_{\ell = 0}^{n-1} (2\ell +
   \alpha + \beta + 2)
 \left\{{N-1 \choose \ell} {(\alpha + \beta + 2)_\ell \over 
   (\alpha + \beta + N + 2)_\ell} {h_\ell (x, \alpha, 
   \beta + 1, N - 1)\over \beta + 1}\right\}^2\cr
& \cr
&+~ (1 - {x\over N}) \left\{{N-1\choose n} {\beta + n + 1\over \beta
   + 1} {(\alpha + \beta + 2)_n \over (\alpha + \beta + N + 2)_n} h_n
   (x, \alpha, \beta + 1, N - 1)\right\}^2\cr
& \cr
&= ~1 ~-~ {x\over N} \left\{{N-1\choose n} {(\alpha + \beta + 2)_n\over 
   (\alpha + \beta + N + 2)_n} h_n (x-1, \alpha + 1, \beta,
   N-1)\right\}^2\ \cr
}
$$
\bigskip
\noindent
{\bf \underbar{Proof.}}~~
We will only give a proof of the identity (a) using the general
result in (2.4).
All other cases are treated similarly where the identity (2.4) has
to be replaced by the corresponding results in Dette (1994).
Observing (2.4), (3.1), (3.2) and (3.3) we have to find the
determinants $\underline{D}_{2\ell} (\xi_\rho),
\overline{D}_{2\ell} (\xi_\rho), \underline{D}_{2\ell - 1}
(\xi_\rho), \overline{D}_{2\ell - 1} (\xi_\rho)$ where
$\xi_\rho$ is the probability measure corresponding to the jump
function (1.2).
But these determinants can easily be calculated from the leading
coefficients of the orthonormal polynomials with respect to the
measures $d\xi_\rho (x), xd \xi_\rho (x), (N-x) d\xi_\rho (x), x
(N-x) d \xi_\rho (x)$ (see e.g.\ Karlin and Studden (1966), p.\
110).
For example, the orthonormal polynomial with respect to the measure
$x (N-x) d\xi_\rho (x)$ is the Hahn polynomial given in (3.1) and
the leading coefficient is  obtained from the definition of
the Hahn polynomials in terms of the hypergeometric series (see
Section 1).
Thus we have for the leading coefficient of the polynomial in (3.1)
$$
\eqalign{
\sqrt{{(\alpha + \beta + 2) (\alpha + \beta + 3)\over N (N-1)
   (\alpha + 1) (\beta + 1)} \pi_n (\alpha + 1, \beta + 1, N-2)}& \cdot
   {(n + \alpha + \beta + 3)_n\over (\alpha + 2)_n (-N + 2)_n}\cr
& =~ (-1)^n \cdot \sqrt{{\overline{D}_{2n} (\xi_\rho)\over
   \overline{D}_{2n+2} (\xi_\rho)}}\cr
}
$$
or equivalently (using (1.3))
$$
{\overline{D}_{2n+2} (\xi_\rho)\over \overline{D}_{2n}
   (\xi_\rho)}~ =~ {n!~(\alpha + 1)_{n+1} (\beta + 1)_{n+1}  (N + \alpha +
   \beta + 2)_n (N-n-1)_{n+2}\over (\alpha + \beta + 2)_{n+1} (n + \alpha
   + \beta + 3)_n (n + \alpha + \beta + 3)_{n+1}}~.
$$
Similarly we obtain for the ratio of $\underline{D}_{2n}
(\xi_\rho)$ and $\underline{D}_{2n-2} (\xi_\rho)$
$$
{\underline{D}_{2n} (\xi_\rho)\over \underline{D}_{2n-2}
   (\xi_\rho)}~ =~ {n!~(\alpha+1)_n (\beta+1)_n (\alpha + \beta + N + 2)_n
    (N-n + 1)_n\over (\alpha + \beta + n + 1)_{n+1} (\alpha + \beta
   + n+1)_n (\alpha + \beta + 2)_{n-1}}
$$
and a straightforward computation yields
$$
{\overline{D}_{2n} (\xi_\rho)\over \underline{D}_{2n}
   (\xi_\rho)}~ =~ {(N-n) (\alpha + \beta + n + 1)\over n (N+ \alpha +
   \beta + n+1)} {\overline{D}_{2n-2} (\xi)\over \underline{D}_{2n-2}
   (\xi)} ~=~ {(N-n)_n (\alpha + \beta + 2)_n\over n! (N + \alpha +
   \beta + 2)_n}~.
\leqno(3.4)
$$
In the same way we find
$$
{\underline{D}_{2n-1} (\xi_\rho)\over \overline{D}_{2n-1}
   (\xi_\rho)}~ =~ {(\alpha + 1)_n\over (\beta + 1)_n}~,\ ~~~ 
   {\underline{D}_{2n} (\xi_\rho)\over \underline{D}_{2n}^+
   (\xi_\rho)}~ =~ {n\over N} {\alpha + \beta + N + n + 1\over \alpha +
   \beta + 2n + 1}
\leqno(3.5)
$$
and
$$
{\overline{D}_{2n} (\xi_\rho)\over \underline{D}_{2n}^+
   (\xi_\rho)}~ =~ {(\alpha + \beta + 2)_n (N-n)_n\over (n-1)! ~(N +
   \alpha + \beta + 2)_{n-1} (\alpha + \beta + 2n + 1) N}
\leqno(3.6)
$$
where we have used the representation (2.3) and (3.4). The orthonormal
polynomials with respect to the measures $(N-x)d\xi_\rho (x)$ and $
x(N-x)d\xi _\rho (x)$ are given by (3.3) and (3.1) and the 
 assertion (a) of Theorem 3.1 follows now from (2.4), (3.4),
(3.5), (3.6) and straightforward  but tedious algebra.
\hfill \bull
\bigskip
\noindent
The Jacobi polynomials $P^{(\alpha, \beta)}_\ell (x)$ orthogonal
with respect to the (continuous) measure $(1-x)^\alpha (1 +
x)^\beta d x$ and with leading coefficient ${ 2^{- \ell}}
{{2 \ell + \alpha + \beta} \choose \ell}$ can be obtained as
limits from the Hahn polynomials
\hfuzz=8pt
$$
P^{(\alpha , \beta)}_n (x) ~=~ \lim_{N \rightarrow \infty}
{{n + \alpha} \choose \alpha} Q_n (N {{1-x} \over 2}, \alpha,
\beta, N)
\leqno(3.7)
$$
\noindent
and  replacing $x$ by $-x$ it
 is straightforward to show that for the limit (3.7)
Theorem 3.1 gives the corresponding formulas for the Jacobi
polynomials in Dette (1994).  For these polynomials it is well
known that
$|P^{(\alpha,\beta)}_n (x)|$ is bounded by $\max \{|P^{(\alpha,
\beta)}_n (-1)|, |P^{(\alpha, \beta)}_n (1)|\}$  $(n \in \N)$
 if $\max \{\alpha, \beta\} > -{1 \over 2}$.
An upper bound but not necessarily sharp bound for arbitrary
parameters is given by Erd\'elyi, Magnus and Nevai (1992).
For the Hahn polynomials the situation is more complicated.
Zaremba (1975) showed that 
$$
|Q_n (x, \alpha, \beta, N)| ~\leq ~1 
\leqno (3.8)
$$
 for $x=0,\ldots ,N$ provided that 
 $ \alpha \geq \beta > -1$, $n(n+1) \leq N $ and
$$
\alpha^2 + \beta^2 - \alpha \beta + \alpha + \beta ~\geq ~0.
\leqno(3.9)
$$
\noindent
In the following theorem we will give an alternative bound for
these polynomials, where the restriction on the degree of the
polynomials satisfying (3.8) depends on
the parameters of the weight function (1.2) and the inequality
holds for all $x \in [0,N]$.
\bigskip
\medskip
\noindent
{\bf\underbar {Theorem 3.2.}}~~Let $\alpha + \beta > -1$ and
$$
 n(\alpha, \beta, N) ~: =~ \ - {1 \over 2} \{ (\alpha + \beta
 -1) \ - \ 
\sqrt{
(\alpha + \beta + 1) (\alpha + \beta + 2N + 1)}
\},
\leqno(3.10)
$$
\noindent
then the $n$th Hahn polynomial satisfies for all $x \in [0,N]$
 and all $n \leq n(\alpha,\beta,N) $ the inequality
$$
|Q_n (x, \alpha, \beta, N)| ~\leq ~\max \ \left\{ 1, 
{{(\beta + 1)_n}
\over
{(\alpha + 1)_n}} \right\}~ = ~\max \ 
\{| Q_n (0, \alpha, \beta, N)|~,~ |Q_n (N, \alpha,\beta, N)|\}
~.
$$
\bigskip
\noindent
{\bf \underbar {Proof:}}~~The second identity
follows from Karlin McGregor (1961), equation (1.13)
and (1.14).
  Let $\beta \geq \alpha$ and $\alpha + \beta
> -1$, by (3.10) all terms on the left hand side of the identity
in Theorem 3.1(c) are positive  which yields 
(here we replace $n$ by  $n-1$
in 3.1c)) 
$$
|Q_n (x, \alpha, \beta, N)| ~\leq ~ 
{{(\beta + 1)_n}
\over
{(\alpha + 1)_n}}
$$
\noindent
for all $x \in [0,N]$.  If $\alpha \geq \beta$ we use the
symmetry relation
$$
Q_n (x, \alpha, \beta, N) ~=~ (-1)^n
{{(\beta + 1)_n}
\over
{(\alpha + 1)_n}}
Q_n (N-x, \beta, \alpha, N)
$$
\noindent
(see e.g. Nikifarov, Suslov and Uvarov (1991), equation (2.4.18),
or Karlin and McGregor (1961), equation (1.15),
but note that both references use a different notation)
 and obtain from
the first part of the proof
$$
|Q_n (x, \alpha, \beta, N)| ~=~ \left|
{{(\beta + 1)_n}
\over
{(\alpha + 1)_n}}
Q_n (N-x, \beta, \alpha, N) \right| ~\leq ~1
$$
\noindent
for all $x \in [0,N]$.  This completes the proof of the theorem.
\hfill \bull
\bigskip
\noindent
{\bf \underbar {Remark 3.3.}}~~Zaremba (1975) proved (3.8)
for $\alpha \geq \beta > -1$ satisfying (3.9),$
 n(n+1) \leq N$ but only for
the integers $x = 0, \ldots, N$ while Theorem 3.2 gives 
 the sup-norm of the Hahn polynomials 
 for all $\alpha + \beta > -1$.  By restricting on the
set $\{0,1,\ldots,N\}$ and $\alpha \geq \beta > -1$ Zaremba's
bound on the degree of the polynomials (such that (3.8) is
satisfied) is comparable with (3.10).  If $\alpha = \beta = 0$ we
obtain from Zaremba (1975) that (3.8) holds for all $n \leq (-1 +
\sqrt {4N+1})/2$ while Theorem 3.2 establishes the (for $
N\geq 13$ weaker) bound
$(1 + \sqrt{2N+1})/2$.  This can be explained by the fact that
Zaremba's approach is directly related to the discrete Legendre
polynomials $Q_n (x, 0, 0, N)$  (and to the integers $\{0,\ldots ,
N\}$) and the general case is obtained
using a projection formula and results of Askey and Gasper (1971)
(for this step the
 condition (3.9) is used).  However, in most cases
Theorem 3.2 will provide a better bound on the degree of the Hahn
polynomials such that  (3.8) is satisfied.
Furthermore the condition (3.9)
is not needed for establishing these bounds.  For example, if
$\alpha + \beta \geq 1$ and $N\ge 3$, then it is easy to see that $(-1 +
\sqrt{4 N+1})/2 \leq n(\alpha, \beta, N)$ and consequently
Theorem 3.2 gives a better bound on the degree of the
polynomials, compared to the results of Zaremba (1975). Moreover if $n
\leq n(\alpha, \beta, N)$, (3.8) is satisfied for all $ x\in
[0,N].$  As further example consider the case
 $\alpha = \beta > -{1 \over 2}$ and $\beta(\beta +
2) < 0$, then (3.9) is not satisfied and Zaremba's results can
not be applied.  However, we obtain readily from Theorem 3.2 that
(3.8) holds for all $x \in [0,N]$ whenever $n \leq \{-(2\beta - 
1) + \sqrt{(2\beta + 1)(2\beta + 1 + 2N)}\}/2.$
\smallskip
\noindent
Zaremba (1975) considered also the example
$$
Q_n (2, -{1 \over 2}, -{1 \over 2}, n^2) = -{5 \over 3} \ \ \ \
\ (n \geq 2)
\leqno(3.11)
$$
\noindent
in order to show that the condition (3.9) can not be relaxed.  In
this case Theorem 3.2 is not applicable and (3.11) indicates that
the Hahn polynomials $Q_n(x,\alpha,\beta,N)$ may not be bounded
by their absolute values at the points $0$ and $N$ if $\alpha +
\beta \leq -1$.  Nevertheless, the following result provides a
bound for these polynomials without a restriction on their
degree.
\bigskip
\medskip
\noindent
{\bf \underbar {Theorem 3.4.}}~~  Let $\alpha + \beta \leq -1$ and $n
\in \{0, \ldots, N-1\}$ then the Hahn polynomials $Q_n(x, \alpha,
\beta, N)$ satisfy for all $x \in [0,N]$ the inequality
$$
|Q_n (x, \alpha, \beta, N)| ~\leq ~\max \left\{~1~,~
{{(\beta+1)_n}
\over
{(\alpha + 1)_n}}
\right\} \cdot
{{(\alpha + \beta + 2 + N)_{n-1}}
\over
{(\alpha + \beta + 2)_{n-1}}}
{{(n-1)!
\over
(N-n+1)_{n-1}}}.
\leqno(3.12)
$$
\bigskip
\noindent
{\bf \underbar {Proof:}}~~Let $\beta \geq \alpha$, then by the
assumptions all terms in the sums of Theorem 3.1(a) are positive.
Consequently we have
$$
\left|{{N-1} \choose {n-1}}
{{(\alpha + \beta + 2)_{n-1}}
\over
{(\alpha + \beta + N + 2)_{n-1}}}
h_n (x, \alpha, \beta, N)\right| ~\leq ~1
$$
\noindent
which is equivalent to (3.12) for $\beta \geq \alpha$.  The case
$\alpha \leq \beta$ is similar as in the proof of Theorem 3.2 and
therefore omitted. \hfill \bull
\bigskip
\noindent
{\bf\underbar {Remark 3.5.}}~~Note that in general the bound (3.12)
can not be improved.  This follows readily from (3.11) for 
$N = 4, n=2$ $  (\alpha = \beta = -{1 \over 2})$ because in this
case  the
right hand side of (3.12) is also given by ${5 \over 3}$. 

\bigskip
\noindent
For 
 $\alpha = \beta = -{1 \over 2}$ we obtain  the discrete analogue of
the Chebyshev polynomials which are  of particular interest and
considered in the following corollary. This result gives a ``discrete''
version of the trigonometric identity (part (a)).
\medskip
\bigskip
\noindent
{\bf \underbar {Corollary 3.6.}}~~Let $T_n (x,N) = Q_n (x, -{1 \over
2}, -{1 \over 2}, N)$ and $U_n (x, N) = Q_n (x, {1 \over 2}, {1
\over 2}, N)$ denote the discrete Chebyshev polynomials of the
first and second kind, respectively, then we have the following
for all $x \in [0,N]$.

\hfuzz=35pt
\noindent
a) For $n=0,\ldots ,N-1$:
$$
\eqalign{
-x(x-{x \over N})& \sum^{n-1}_{\ell=0} (\ell + 1) \left\{
{{4(\ell + 1)} 
\over
{N-1}}
U_l (x-1, N-2)\right\}^2 \cr
&+~T_{n+1}^2(x,N)~+ ~ 
{x \over N} (1 - {x \over N}) \left\{
{{2(n+1)(N-n-1)}
\over
{N-1}}
U_n(x-1, N-2)\right\}^2
~=~ 1~.\cr 
}
$$
\noindent
b) For $n=0,\ldots ,N-1$:
$$
|T_n(x,N)| ~\leq ~\prod^{n-1}_{j=1} \left(1 + {n \over {N-n+j}}\right) \  .
$$
\noindent
c) 
For $0 \leq n \leq  \sqrt{N+1} $~:
$$
|U_n(x,N)| ~\leq ~1 \   . 
$$
\bigskip
\noindent
{\bf \underbar {Remark 3.7.}}~~Observing that the Jacobi polynomials
can be obtained as the limit (3.7) from the Hahn polynomials and
using formula (4.17) in Szeg\"o (1975) it is easy to see that the
first part of Corollary 3.5 yields $(N \rightarrow \infty \ \
x={N \over 2} (1-z))$ the trigonometric identity $(1-z^2)
U^2_n(z) + T^2_{n+1}(z) = 1$ for the Chebyshev polynomial of the
first and second kind while part (b) and (c) establish the bounds
$|T_n(z)| \leq 1, |U_n(z)| \leq n+1$ $ (z \in [-1,1])$ for these
polynomials (note that $\lim_{N\to \infty} U_n({N\over 2}(1-z),N)
={U_n(z)\over n+1}$).

\bigskip
\noindent
We will conclude this section with a brief discussion of related 
 results for the Hahn- Eberlein and the dual Hahn
polynomials. The Hahn- Eberlein polynomials
 are obtained from the Hahn polynomials $Q_n(x,
\alpha, \beta, N)$ for $\alpha < -N, \beta < -N$ (see e.g. Rahman
(1978) or Eberlein (1964)).  For this choice the mass function in
(1.2) still defines a probability measure on $\{0,\ldots,N\}$
and consequently the orthogonal polynomials
$Q_n(x,\alpha,\beta,N)$ with respect to this measure are well
defined and called Hahn- Eberlein polynomials.  These polynomials
have some applications in coding theory (see e.g. Sloane (1975)).
Obviously, the identities of Theorem 3.1 can be extended to the
region $\alpha < -N, \beta < -N$ and as a consequence we obtain
the following bound for the Hahn- Eberlein polynomials.
\bigskip
\noindent
{\bf \underbar {Theorem 3.8.}}~~Let $\alpha < -N, \beta < -N, \alpha
+ \beta < -2N-1$ and
 $$
\tilde n (\alpha, \beta, N) \ = \ -{1 \over 2}
\{ (\alpha + \beta - 1) \ + \
\sqrt{(\alpha+\beta+1)(\alpha+\beta+1+2N)}\} .
$$
\noindent
The Hahn- Eberlein polynomials $Q_n(x,\alpha,\beta,N)$ satisfy for
all $x \in [0,N]$ and $n \leq $\~n$ (\alpha, \beta, N)$ the
inequality
$$
|Q_n(x,\alpha,\beta,N)| ~\leq ~\max \left\{1, 
{{(\beta+1)_n} \over {(\alpha+1)_n}}\right\} = \max
\{|Q_n(0,\alpha,\beta,N)|, |Q_n(N,\alpha,\beta,N)|\} ~.
$$
\bigskip
\noindent
The dual Hahn polynomials $R_k(x,\alpha,\beta,N)$ ($\alpha,\beta > -1$)
 are related to the Hahn polynomials
by the equation
$$
R_k(x(x+\alpha+\beta+1))~=~ Q_x(k,\alpha,\beta,N)
$$
($k,x=0,\ldots ,N$) and are orthogonal
on the interval $[0,N(N+\alpha+\beta+1)]$.
 For a detailed discription of these polynomials
including the recurrence relation and the orthogonality relation
	we refer the reader to the work of
Karlin and McGregor (1961). By a similar analysis as in the proof 
of Theorem 3.2
we obtain the following bound for these polynomials.	
\bigskip
\noindent
{\bf \underbar{Theorem 3.9.}}~~Let $\alpha,\beta > -1$,  $N+1+
\beta -\alpha \geq 0$ and 
$$
n^*(\alpha,\beta,N)~=~{1\over 2}\min\{N+2,~ N+1+\beta -\alpha\}
.
$$
If $n \leq n^*(\alpha,\beta,N)$ then the dual Hahn polynomial 
$R_n(x,\alpha,\beta,N)$ satisfies for all $x \in [0,N(N+\alpha
+\beta+1)]$ the inequality
$$
|R_n(x,\alpha,\beta,N)|~\leq ~ {(N+1+\beta-n)_n \over
(\alpha+1)_n}~=~|R_n(N(N+\alpha+\beta+1),\alpha,\beta,N)|~~.
$$
\bigskip
\noindent
{\bf \underbar{Proof:}}~~Let $\xi_D$ denote the measure which puts
masses
$$
\xi_D(\lambda_x) ~=~{\pi_x(\alpha,\beta,N) \rho (0)}
$$ 
at the points $\lambda_x = x(x+\alpha+\beta+1)$ ($x=0,\ldots ,N$)
where $\pi_x(\alpha,\beta,N)$ and $\rho(x)=\rho (x,\alpha,
\beta,N)$ are defined in (1.3)
and (1.2), respectively. By the results of Karlin and McGregor
 (1961) (equation (1.20)) it follows that $\xi_D$ defines a 
probability measure on the interval $[0,N(N+\alpha+\beta+1)]$
	and that the orthonormal polynomials with respect 
to $d\xi_D(x) $ are given by
$$
P_l(x)~=~
\sqrt{\rho (l) \over \rho (0)} R_l(x,\alpha,\beta,N)
~~~~~~~~~(l=0,\ldots ,N). \leqno (3.13)
$$ 
According to Theorem 3.1 in Dette (1994) it follows that  the orthonormal
polynomials $P_l(x)$, $Q_l(x)$ and $S_l(x)$ with respect
to
the measures $d\xi_D(x)$, $x[N(N+\alpha+\beta+1)-x]d\xi_D(x)$
and $[N(N+\alpha +\beta +1)-x]d\xi_D(x)$ satisfy the identity
$$
\eqalign{
& x[N(N+\alpha +\beta +1)-x] \sum^{n-1}_{\ell=0}
{{{\underline D}_{2 \ell + 1} (\xi_D)}
\over
{{\overline D}_{2 \ell + 1} (\xi_D)}}
\left[
{{{\underline D}_{2 \ell } (\xi_D)}
\over
{{\overline D}_{2 \ell } (\xi_D)}}
-
{{{\underline D}_{2 \ell+2} (\xi_D)}
\over
{{\overline D}_{2 \ell+2} (\xi_D)}}
\right]
Q^2_\ell (x) \cr
&~~~~~~+~x[N(N+\alpha +\beta +1) -x]
{{{\underline D}_{2n} (\xi_D) \ {\underline D}_{2n+1} (\xi_D) \
{\overline D}_{2n+2} (\xi_D)}
\over
{{\overline D}_{2n} (\xi_D) \ {\overline D}_{2n+1} (\xi_D) \
{\overline D}^-_{2n+2} (\xi_D)}}
Q^2_n (x) \cr
&~~~~~~+ ~[N(N+\alpha +\beta +1)-x] \sum^{n}_{\ell=0}
{{{\underline D}_{2 \ell} (\xi_D)}
\over
{{\overline D}_{2 \ell} (\xi_D)}}
\left[
{{{\underline D}_{2 \ell -1} (\xi_D)}
\over
{{\overline D}_{2 \ell -1} (\xi_D)}}
- 
{{{\underline D}_{2 \ell +1} (\xi_D)}
\over
{{\overline D}_{2 \ell +1} (\xi_D)}}
\right]
S^2_\ell (x) \cr
& ~~~~~~ 
=~
1 ~-~ 
{{{\underline D}_{2n+1} (\xi_D) {\underline D}_{2n+2} (\xi_D)}
\over
{{\overline D}_{2n+1} (\xi_D) {\overline D}_{2n+2}^- (\xi_D)}}
P^2_{n+1} (x)
}
\leqno (3.14)
$$
($n=0,\ldots ,N-1$). By a similar reasoning as in the proof of Theorem 3.1
we obtain for the ratios of the determinants in (3.14)
$$
{{\underline D}_{2l} (\xi_D)\over {\overline D}_{2l}(\xi_D)}~-~
{{\underline D}_{2l+2}(\xi_D)\over {\overline D}_{2l+2}(\xi_D)}~=~
{l!\over (N-l-1)_{l+1}}(N-2l-2) ~~~~~~~~~~(l=0,\ldots ,n-1),
$$
$$
{{\underline D}_{2l-1} (\xi_D)\over {\overline D}_{2l-1}(\xi_D)}~-~
{{\underline D}_{2l+1}(\xi_D) \over {\overline D}_{2l+1}(\xi_D)}~=~
{(\alpha +1)_l \over (N+\beta -l)_{l+1} }
(N-1+\beta -\alpha -2l) ~~~~~~~~~~~~~~~~(l=0,\ldots ,n)
$$
and
$$
\eqalign{
{{\underline D}_{2n+1}(\xi_D) { \underline D}_{2n+2}(\xi_D) \over
{\overline D}_{2n+1}(\xi_D) {\overline D}_{2n+2}^- (\xi_D)}P_{n+1}^2(x)
~&=~
{(\alpha +1)_{n+1} (n+1)! \over 
(N+\beta -n)_{n+1} (N-n)_{n+1}} {\rho (n+1) \over \rho (0)}
R_{n+1}^2 (x,\alpha,\beta,N) \cr
&=~\left( {(\alpha +1)_{n+1} \over
(N+\beta -n)_{n+1}} R_{n+1}(x,\alpha,\beta,N) \right)^2
\cr
}
$$
where we have used (3.13) and (1.2)
in the last identity. By the assumptions of the theorem
 all terms on the left hand side in (3.14) are positive and the assertion
follows from
$$
R_n(N(N+\alpha+\beta+1))~=~Q_N(n,\alpha,\beta,N)~=~ (-1)^n {(N+1+\beta -n)_n 
\over (\alpha +1)_n} 
$$
which can easily be proved by an induction argument.\hfill \bull
\bigskip
\bigskip
\noindent
{\bf \underbar {4. \ Krawtchouk polynomials.}}~~In this Section we
will apply the results of Section 2 and 3 in order to obtain
similar results for the Krawtchouk polynomials.  We will mainly
use the representation (1.4) of $k_n(x,p,N)$ as the limit of the
Hahn polynomials $Q_n(x,\alpha,\beta,N)$ when $\alpha = pt, \beta
= qt$ and $t \rightarrow \infty$.  By this relation the following
results are immediate consequences of Theorem 3.1 and 3.2.
\bigskip
\noindent
{\bf\underbar {Theorem 4.1.}}~~For $\ell = 0, \ldots, N$ define
$\tilde k_\ell (x,p,N) = {N \choose \ell} ({p \over q})^\ell k_\ell
(x,p,N)$.  The Krawtchouk polynomials satisfy the following
identities.

\noindent
a) For $n=0,\ldots ,N-1:$
$$
\eqalign{
\sum^{n-1}_{\ell=1} ({{2 \ell} \over N} - 1) \{ \tilde k_\ell
(x,p,N)\}^2 &\ + \ \{ {n \over N} \tilde k_n (x, p, N)\}^2
~+~ (1 - {x \over N}) 
{{q-p} 
\over 
{q^2}} 
\sum^{n-1}_{\ell=0}
\{\tilde k_\ell (x,p,N-1)\}^2
\cr 
&~=~ 1 ~-~ {x \over N} (1 - {x \over N}) \left\{ {{k_{n-1}(x-1,p,N-2)} \over
q} \right\}^2 \  . \cr
}
$$
\noindent
b) For $n=0,\ldots ,N-1$:
$$
\eqalign{
& \sum^n_{\ell=1} ({{2 \ell} \over N}-1) \left\{ {N \choose \ell} k_\ell
(x, p, N) \right\}^2 \ + \ {x \over N} \left\{ {{N-1} \choose n} k_n (x-1,
p, N-1) \right\}^2
\cr & \cr
&+~ {{p-q} \over {p^2}} {x \over N} \sum^{n-1}_{\ell=0} \left\{ {{N-1}
\choose \ell} k^2_\ell (x-1, p, N-1) \right\}^2 \cr
&\cr
&
~=~ 1~ -~ (1 - {x \over N}) 
\left\{ {{N-1} \choose n} k_n (x, p, N-1) \right\}^2 ~.\cr
} 
$$
\smallskip
\noindent
c) 
For $n=0,\ldots ,N-2$:
$$
\eqalign{
&x (1-{x \over N}) \sum^{n-1}_{\ell=0} (N-2\ell-2) \left\{ 
{{p^\ell}
\over
{q^{\ell+1}}}
{{k_\ell(x-1, p, N-2)}
\over
{(N-1)}}
\right\}^2 \cr 
& \cr 
&~+~ {{q-p}
\over
{q^2}} (1 - {x \over N}) \sum^n_{\ell=0} 
\left\{ {{p^\ell} \over {q^\ell}} 
k_\ell (x,p,N-1) \right\}^2 
\cr & \cr
&+~ {x \over N} (1 - {x \over N}) \left\{
{{(N-n-1)p^n k_n (x-1,p,N-2)}
\over
{q^{n+1} (N-1)}}
\right\}^2
~=~ 1 ~- ~\left\{ 
{{p^{n+1}}
\over
{q^{n+1}}} 
k_{n+1} (x, p, N) \right\}^2 
}
$$
\smallskip
\noindent
d) For $n=0,\ldots ,N-1$:
$$
\eqalign{
&\sum^n_{\ell=1} ({{2 \ell} \over N} - 1) \{\tilde k_\ell (x, p,
N) \}^2 +\
{{q-p} \over {q^2}} (1 - {x \over N}) \sum^{n-1}_{\ell=0}
\{\tilde k_\ell (x,p,N-1)\}^2
\cr &\cr 
&~+~ (1-{x \over N}) \{ \tilde k_n (x, p, N-1)\}^2 ~=~ 1 ~-~ {x \over N}
\{\tilde k_n(x-1, p, N-1)\}^2 \  .\cr
} 
$$

\bigskip
\bigskip
\noindent
{\bf \underbar {Theorem 4.2.}}~~Let $n \leq {N \over 2}+1$, then the
$n$th Krawtchouk polynomial $k_n(x,p,N)$ satisfy for all $x \in [0,N]$ the
inequality
$$
|k_n (x, p, N)| ~\leq ~\max \left\{ 1, ({q \over p})^n \right\} ~=~ \max 
\{| k_n (0, p, N)|~, ~ |k_n (N, p,N)|\}~.
$$
\bigskip
\bigskip
\noindent
{\bf Acknowledgements.}~~Parts of this paper were written while the author was
visiting the University of G\"ottingen and Purdue University, West
Lafayette. The author
would like
to thank the Institut f\"ur Mathematische Stochastik and the
 Department of Statistics for the hospitality
 and the Deutsche Forschungsgemeinschaft
for the financial support that made a visit to the United States
possible. I am also indebted to Dick Askey and George Gasper
for their helpful comments and help with the references.
It was George Gasper who convinced me that identities of the
form (2.4) may give useful bounds for orthogonal polynomials
while Dick Askey communicated to me the results of S.K. Zaremba.

\bigskip
\bigskip
\noindent
{\bf REFERENCES}
{\advance\leftskip by 0.2in \parindent=-0.2in

R. Askey, G. Gasper (1971).  Jacobi polynomial expansions of
Jacobi polynomials with nonnegative coefficients, Proc. Camb.
Phil. Soc., 70, 245-255.

H. Dette (1994).  New identities for orthogonal polynomials on
compact intervals, J. Math. Anal. Appl., to appear.

P. J. Eberlein (1964).  A two parametric test matrix.  Math.
Comput., 18, 296-298.

T. Erd\'elyi, A.P. Magnus, P. Nevai (1992). Generalized Jacobi weights,
Christoffel functions, and Jacobi Polynomials. Ohio State Mathematical Research
Institute Preprints, \#92-29.

G. Gasper (1974).  Projection formulas for orthogonal polynomials
of a discrete variable.  J. Math. Anal. Appl., 45, 176-198.

G. Gasper (1975).  Positivity and special functions.  Theory and
Application of Special Functions (ed. R. Askey), Academic Press
Inc., 375-433.

W. Hahn.  \"Uber Orthogonalpolynome, die q -
Differenzengleichungen gen\"ugen, Math. Nachrichten, 2, 4-34.

S. Karlin, J. McGregor (1961).  The Hahn polynomials, formulas
and an application.  Scripta Math. 26, 33-46.

S. Karlin, J. McGregor (1962).  On a genetics model of Moran,
Proc. Cambridge Phil. Soc., 58, 299-311.

S. Karlin, L. S. Shapeley (1953).  Geometry of Moment Spaces,
Amer. Math. Soc. Memoir No 12, Amer. Math. Soc., Providence.

S. Karlin, W. J. Studden (1966).  Tchebycheff Systems:  with
applications in analysis and Statistics, Interscience Publ., NY.

M. Krawtchouk (1929).  Sur une g\'en\'eralisation des polynom\'es
d' Hermite.  Comptes Rendus de l' Acad\'emie des Sciences, Paris,
189, 620-622.

A. F. Nikiforov, S. K. Suslov, V. B. Uvarov (1991).  Classical
Orthogonal Polynomials of a Discrete Variable, Springer Verlag,
New York.

M. Rahman (1978).  A positive kernel for Hahn-Eberlein
polynomials. SIAM J. Math. Anal., 9, 891-905.

N. J. A. Sloane (1975).  An introduction to association schemes
and coding theory.  Theory and Application of Special Functions
(ed. R. Askey),  Academic Press, New York, 225-260.

G. Szeg\"o (1975).  Orthogonal Polynomials, Amer. Math. Soc.
Colloq. Publ., Vol. 23, Amer. Math. Soc., NY.

M. W. Wilson (1970).  On the Hahn polynomials, SIAM J. Math.
Anal., 1, 131-139.

S. K. Zaremba (1975).  Some properties of polynomials orthogonal
over the set $<1,2,\ldots,N>$, Ann. Mat. Pura Appl., 105,
333-345.

} 

\bye